\documentclass[11pt]{article}
\usepackage{amssymb,amsmath,mathrsfs,amsfonts,graphicx,epsf} 
\usepackage{color}
\usepackage{xcolor}

 \usepackage{fullpage}

\usepackage{authblk}

\newtheorem{theorem}{Theorem}
\newtheorem{corollary}{Corollary}
\newtheorem{lemma}{Lemma}
\newtheorem{proposition}{Proposition}
\newtheorem{definition}{Definition}
\newtheorem{remark}{Remark}

\newcommand{\bt}{\begin{theorem}}
\newcommand{\et}{\end{theorem}}
\newcommand{\bl}{\begin{lemma}}
\newcommand{\el}{\end{lemma}}
\newcommand{\bp}{\begin{proposition}}
\newcommand{\ep}{\end{proposition}}
\newcommand{\bc}{\begin{corollary}}
\newcommand{\ec}{\end{corollary}}
\newcommand{\bdeff}{\begin{definition}}
\newcommand{\edeff}{\end{definition}}
\newcommand{\brem}{\begin{remark}}
\newcommand{\erem}{\end{remark}}
\newcommand{\bqnn}{\begin{equation*}}
\newcommand{\eqnn}{\end{equation*}}

\renewcommand{\r}[1]{(\ref{#1})}

\newcommand{\bi}{\begin{itemize}}
\newcommand{\iii}{\item}
\newcommand{\ei}{\end{itemize}}
\newcommand{\bpm}{\begin{pmatrix}}
\newcommand{\epm}{\end{pmatrix}}
\newcommand{\bd}{\begin{description}}
\newcommand{\ed}{\end{description}}

\newcommand{\bqn}{\begin{eqnarray}}
\newcommand{\eqn}{\end{eqnarray}}

\newcommand{\ba}[1]{\begin{array}{#1}}
\newcommand{\ea}{\end{array}}

\newcommand{\R}{\mathbb{R}}

\newcommand{\lam}{\lambda}

\newcommand{\g}{\gamma}

\newcommand{\ph}{\varphi}

\newcommand{\con}{{\mathcal C}}

\newcommand{\la}{\left\langle}
\newcommand{\ra}{\right\rangle}

\newcommand{\sign}{\mathrm{sign}}

\title{\LARGE \bf
Optimal control of a multi-scale dynamic  model  for biofuel production
}

\author{
Roberta~Ghezzi\footnote{
Institut de Math\'ematiques de Bourgogne, UBFC, 9 Avenue Alain Savary, 21078 Dijon, France {\tt roberta.ghezzi@u-bourgogne.fr}}
\quad
Benedetto Piccoli\footnote{Department of Mathematical Sciences and Center for Computational and Integrative Biology, Rutgers University 311 N 5$^{\textrm{th}}$ St,   08102 Camden NJ, USA {\tt piccoli@camden.rutgers.edu}}
}

\begin{document}

\maketitle

\begin{abstract}
Dynamic flux balance analysis of a bioreactor is based on the coupling between a dynamic problem, which models the evolution of biomass, feeding substrates and metabolites, and a linear program, which encodes the metabolic activity inside cells. We cast the problem in the language of optimal control and  propose a hybrid formulation  to model  the full coupling between  macroscopic and microscopic level. On a given location of the hybrid system we analyze necessary conditions given by the Pontryagin Maximum Principle and discuss the presence of singular arcs. In particular, for the single-input case we prove that optimal controls are bang-bang. For the multi-input case, under suitable assumptions, we prove that generically with respect to initial conditions optimal controls are bang-bang.
\end{abstract}

\section{Introduction}

Biofuels provide a concrete  answer to the pressing need for renewable energy. The problem of increasing the efficiency and reducing the cost of biofuel production  has  been subject of intensive research \cite{optcontbiofi,optcontbiof}. A feasible source to obtain biofuels consists in using ethanol produced by  cyanobacteria and microalgae \cite{hensoneth,hensonethanol}.   This paper deals with the application of optimal control techniques to a general bioreactor for biofuel production. 
  
Since the last decades, optimization of bioprocesses has been a line of research connecting optimal control to system biology,  see \cite{bioprocess-review} and references therein. 
 Progress in plant genetic engineering has opened novel opportunities to use plants as bioreactors for safe and cost effective production of vaccine antigens. A review of
methods and applications of plant, tissue and cell culture    based expression strategies and their use as bioreactors for large scale production of pharmaceutically important proteins can be found in \cite{vaccine}.
Exploiting bioreactors was also proved a fruitful method to deal with the problem of  wastewater treatment in environmental engineering \cite{moreno}. 
In this context  a   dynamic optimization problem arising frequently  is the minimization of the time needed to reach a fixed target configuration for the bioreactor, see \cite{Rapaport}. In this case,  the trajectory evolves according to a system of nonlinear ODEs and may satisfy  some boundary constraints  as well.  The right-hand side of the evolution equations involves not only control functions, i.e.  parameters through which we can  modify the dynamics, but also  some unknown functions (e.g. biomass growth rate) of the evolving quantities themselves.    Therefore, computing optimal controls through standard techniques (necessary conditions) is affected by the unknown functions' behavior. An interesting approach to deal with this issue using observability techniques was proposed in \cite{gauthier-busvelle}, where an application to bioreactors is also provided (see also \cite{gauthier-bior}).

Other applications of optimal control techniques have been studied in the field of biological systems
 \cite{Alford}. In the context of medical treatment,    chemotherapy was used \cite{hiv} as a dynamic control  to optimize treatment scheduling in early stage HIV-infected cases. In that paper the authors use the effect of chemotherapy on viral production to maximize   benefits in terms of T cell count and minimize the systemic cost of the treatement.    
  Time dependent control strategies were also exploited in models to contain the emergence of drug-resistant strains of tubercolosis  \cite{tb}. Here controls are represented by efforts in finding patients in which virus is only latent and in completing treatment for patients in which virus is already active. The objective function  balances   the effect of minimizing the cases of latent and infectious drug-resistant tubercolosis and minimizing the cost of implementing the control treatments.

In this paper we are concerned with the optimization of a bioprocess for biofuel production.   Namely, we consider a model for the metabolic activity of  a microorganism (e.g. {\it E. Coli} or {\it Saccharomyces cerevisiae}) in a fed-batch culture with different feeding substrates. The optimal control problem is to maximize  the productivity of a  certain side metabolite (e.g., ethanol)  through different feeding rates.

A steady-state approach to model cellular metabolism is Flux Balance Analysis (FBA, see \cite{palsson-book}). The main assumption of FBA is that metabolic activity of cells is performed in such a way that the growth rate of cells is maximized.   Since  genome-scale stoichiometric models for bacteria such as {\it E. Coli} are available, this translates in a linear program where the objective is cellular growth rate and constraints are given by metabolic reactions. Optimization is then performed by means of genetic manipulations on the bacteria (gene deletions and insertions). 
Other    approaches based on flux balance analysis were proposed that take account   of transcriptional and regulatory effects in \cite{rfba}, that couple the steady-state metabolic activity with a dynamic model \cite{odemodel,dfba} and that integrate both aspects \cite{ifba}.   Optimizing ethanol productivity   in fed-batch cultures modeled through Dynamic Flux Balance Analysis \cite{biofuels} consists in using outputs of FBA (cellular growth rate and metabolite fluxes) to update at each time step the dynamics for evolving extracellular quantities. Within this framework, control can be performed at two levels: intracellular controls (genetic modifications), which are implemented by acting on constraints of the FBA, and extracelullar controls (of dynamic nature), which are implemented through time-dependent feeding rates. 
 Genetic strategies were deeply investigated in \cite{biofuels} for {\it in silico} evolution of a yeast strain in glucose and xylose media to maximize ethanol productivity  with constant feeding rates.
 Besides genetic strategies, another parameter which is used as a control in \cite{hensonethanol,biofuels} is the switching time of oxygen concentration. More precisely, the authors treat dissolved  oxygen concentration as an independent variable, assuming it could be regulated by a feedback controller, and switch from aerobic  to anaerobic growth to promote ethanol production during later stages of the batch. Then, they analyze sensitivity of ethanol productivity (for different genetic strains) to the aerobic-anaerobic switching time, while feeding substrates' rates are kept constant.

\begin{figure}[h!]\label{schema}
\begin{center}
\input{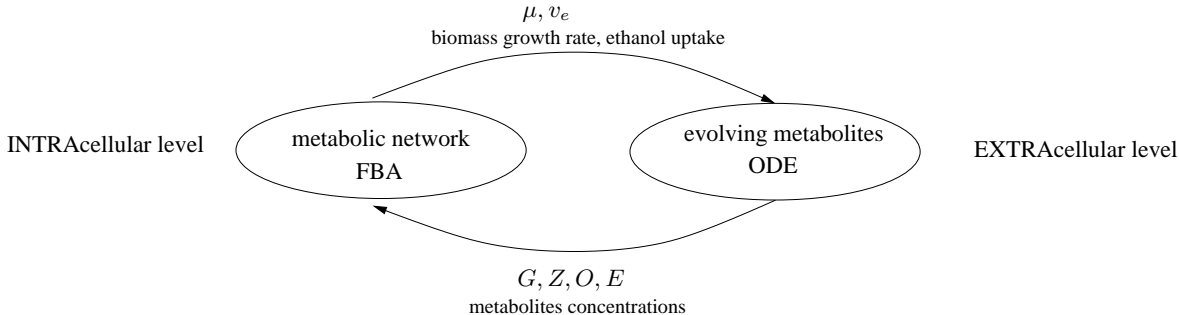}

\caption{Bioprocess scheme exhibiting full coupling between metabolic activity and external dynamics}
\end{center}
\end{figure}

In this paper, we focus on a unified model dealing with   extracellular controls, i.e., feeding rates for glucose, xylose and oxygen, which we allow to be  time dependent. Namely, we cast the problem of optimizing ethanol for the bioreactor in the language of optimal control, where concentrations of feeding substrates are modified by an external agent in order to enhance ethanol production. Even though we do not consider genetic manipulations, our   model  can  integrate intracellular controls. 
Our  purpose is to develop  a general  model featuring a full coupling between extracellular dynamics and intracellular metabolic activity: at each time instant solutions of ODEs provide fluxes to constraint the FBA, whereas FBA gives cellular growth rate and ethanol uptake appearing in the ODEs see Figure~\ref{schema}. 
Since both the objective function and the constraints in FBA are linear (with respect to fluxes), outputs are piecewise linear. Dependence of fluxes on metabolite concentrations is modeled through Michaelis--Menten behavior, i.e., through   rational functions. Therefore, FBA outputs appearing in dynamics are piecewise smooth functions of evolving quantities. The main idea is to implement the coupling between the intracellular and extracellular level through a hybrid control system where in each location the outputs of FBA are smooth.  One of the advantages of a hybrid formulation is that it allows to take account of  different timescales (for the general theory of hybrid control systems we refer the reader to \cite{hyb,hyblibro}). One can either assume that information translates instantaneously from the micro level to the macro one (and viceversa) or one can implement  delays, which are observed in experimental data. This is useful as there are configurations of the extracellular environment (such as saturation of glucose or oxygen) which may cause
a major change in the metabolic pathway involved inside the cell with a consequent delay in the variation of outputs of the FBA.

 As a first step toward a unified model,   we   focus on the analysis of optimal trajectories contained in a region where  FBA ouptuts   are smooth. In other words, we study what happens in each location of the hybrid model.
More precisely,  the location is characterized by  a control-affine system of the type
\bqnn
\dot x= F_0(x,y^x)+u_1 F_1(x)+u_2 F_2(x) +u_3F_3(x), 
\eqnn
where $x$ is a vector containing substrate concentrations (biomass, glucose, xylose, ethanol, oxygen), $u_1,u_2, u_3$ are feed rates of glucose, xylose and oxygen, $y^x$ are parameters coming from FBA, and $F_i$, $i=0,\dots, 3$ are smooth vector fields. We consider an optimal control problem in Mayer form where we optimize the total amount of ethanol at the final time. Since the system is control-affine, the maximization condition of the  Pontryagin Maximum Principle allows to compute 
the control along an extremal trajectory as a feedback law, as long as the corresponding switching function (derivative of Hamiltonian with respect to control) has only isolated zeros. Otherwise, when the trajectory shows a singular arc, i.e., a time interval where a switching function vanishes identically,   the main tool to find controls is exploiting conditions given by  annihilation of higher order time-derivatives. 
 As our interest is driven by applications, one would rather avoid singular arcs, which represent an obstacle to efficient and reliable numerical simulations. 
 
 First, we consider the single-input case, where optimization is performed under the action of one control. In this connection, we prove that every optimal trajectory is a concatenation of arcs where the feeding rate is either absent or maximal. Then we study the general three-input case where ethanol productivity is maximized using glucose, xylose and oxygen feeding rates as controls. The situation is more intricate and some specific singular arcs may indeed be optimal. Under an additional assumption, which amounts to disregard the presence of a preferred substrate, we show that at most initial conditions optimal trajectories do not have  singular arcs and are concatenation of arcs where all feeding rates are either absent or maximal.


 The structure of the paper is the following. In Section~\ref{sec2} we recall the general model of a bioreactor and formulate the optimal control problem where the cost to be maximized is ethanol productivity of the bioreactor. In Section~\ref{secmainth} we recall a classical first order necessary condition for optimal control problems. Then, we consider in Section~\ref{sinput} the single-input case and show that every optimal trajectory is bang-bang.
 In Section~\ref{minput} we deal with the general case of three feeding substrates and analyze possible singular arcs. Finally, in Section~\ref{conclu} we compare our results with the ideas in \cite{biofuels}.

\section{Problem formulation}\label{sec2}

Bioreactors  are processes  where a living microorganism metabolizes some substrates, consequently grows and produces other metabolites. In \cite{biofuels} the authors consider {\it in silico} evolution of a yeast strain, {\it Saccharomyces cerevisiae}, which grows in a fed-batch culture with   glucose and xylose and produces ethanol. 
The dynamic model can be applied to a general bioreactor.

 We denote by  $V$ the total  culture volume, which is assumed to grow linearly with respect to time with  constant rate $F$.
Biomass concentration is denoted by $X$. The total biomass  in the culture evolves  linearly with a growth rate $\mu$ depending on substrates concentrations.
Concentrations of feeding substrates, glucose $G$, xylose $Z$, and oxygen $O$, are characterized by an evolution which takes account of 
a feeding rate and a compensation term  due to metabolism of the microorganism. Feeding rates for substrates represent the control we perform on the system and are denoted by $u_1, u_2, u_3$. The organisms metabolize  glucose, xylose and oxygen  with specific rates (or fluxes) $v_g, v_z, v_o$   depending on substrates concentrations. The produced metabolite under study is ethanol, with concentration $E$.   We assume that the organism produces ethanol proportionally to the total biomass through a specific rate $v_e$.
 Therefore, the control system we analyze is 
\begin{equation}\label{eq:cs}
\left\{\ba{ccc}
\dot{VG}&=&u_1-v_g(G,E)VX,\\
\dot{VZ}&=&Fu_2-v_z(G,Z,E)VX,\\
\dot{VO}&=&Fu_3-v_o(O)VX,\\
\dot{VE}&=&v_e(G,Z,O,E) VX,\\
\dot{VX}&=&\mu(G,Z,O, E) VX,\\
\dot V&=&F.
\ea\right.
\end{equation} 
First of all, 
 oxygen uptake kinetics follows Michaelis--Menten law
 \bqnn
 v_o(O)=v_{o\max}\frac{O}{k_o+O},
 \eqnn
 whereas 
 glucose    uptake kinetics has an additional regulatory term to capture growth rate suppression due to high ethanol concentration, i.e.,
\begin{eqnarray*}
v_g(G,E)&=&v_{g\max}\frac{G}{k_g+G}\frac{1}{1+E/k^g_{ie}}.
\end{eqnarray*}
Xylose uptake kinetics has a similar form with another regulatory term to account for inhibited xylose metabolism in presence of the preferred substrate (glucose),
\bqnn
v_z(G,Z,E)=v_{z \max}\frac{Z}{k_z+Z}\frac{1}{1+E/k^z_{ie}}\frac{1}{1+G/k_{ig}}.
\eqnn
Parameters $v_{o\max},v_{g\max},v_{z\max}, k_o, k_g, k^g_{ie}, k_z,k^z_{ie}$ are positive and constant. 
As for $\mu, v_e$, the model is based on the principle that the metabolic activity of the microorganism is performed so that biomass growth rate is maximized.
This is equivalent to say that   $\mu$ and $v_e$ are outputs of the optimization problem stated as
  \begin{eqnarray}
&&\mu(G,Z,O,E)=\max_{\bar v\in\R^{n}} \,\sum_{j=1}^{n} w_j \bar v_j\label{lp}\\
&&\quad\qquad\qquad\qquad \mbox{ s.t. }  \, S\bar v = 0\nonumber\\
 &&\quad\qquad\qquad\qquad   \, \bar v_g=v_g(G,E) \nonumber\\
 &&\quad\qquad\qquad\qquad   \, \bar v_z=v_z(G,Z,E)\nonumber\\
&&\quad\qquad\qquad\qquad   \, \bar v_o=v_o(O)\nonumber\\
&&\quad\qquad\qquad\qquad \, 0\leq  \bar v_j\leq \tilde v_j \nonumber\\
&& \nonumber\\
 && v_e(G,Z,O,E) =\max\{\bar v_e\mid \bar v\in\mathrm{argmax }\, \mu(G,Z,O,E)\}.\label{ve}
 \end{eqnarray}
 
In \r{lp}, $\bar v\in\R^n$ is the vector of fluxes (among which  glucose, xylose, oxygen and ethanol fluxes) considered in the model; $w\in[0,1]^n$ is the vector of weights which determines fluxes producing biomass;  $S\in M^{r\times n}(\R)$ is the stoichiometric matrix which encodes the metabolic network inside the cell (involving $r$ reactions and $n$ metabolites);   $\tilde v$ is an upper bound associated with the microorganism,   
(see \cite{palsson-book} for an exhaustive treatment of metabolic networks in systems biology.)  For our purposes, $w,S,\tilde v$ are considered as given parameters. They depend  on the specific strain of microorganism used in the bioreactor and are usually determined through experimental data,  see for instance supplementary data of  \cite{ifba}  for {\it E. Coli} or \cite{hensonethanol} for {\it Saccharomyces cerevisiae}.

 Equation \r{ve} is to be read as follows.  The set  $\mathrm{argmax } \, \mu(G,Z,O,E)$ is the set of vectors $\bar v$ that are solutions to \r{lp}, i.e., that satisfy  the constraints   and realize the maximum,  
$\mu(G,Z,O,E)=\sum_j w_j\bar v_j. 
$
Then  $v_e(G,Z,O,E)$ is the maximum of $\bar v_e$, the component of $\bar v$ corresponding to ethanol flux, as $\bar v$ varies in $\mathrm{argmax } \, \mu(G,Z,O,E)$. 

We fix the final time $t_f>0$  and an initial condition $(G_0,Z_0,O_0,E_0,X_0,V_0)$ and we  seek to optimize ethanol productivity along a trajectory of \r{eq:cs}. More precisely, we optimize the cost
\begin{equation}\label{produtt}
\displaystyle{\max_{(u_1,u_2, u_3)\in{\cal U}} \frac{V(t_f)E(t_f)}{\int_0^{t_f}(u_1(s)+u_2(s)+u_3(s))\,ds}}
\end{equation}
where the class of admissible controls is   
\begin{equation*}
{\cal U}=\{(u_1,u_2,u_3):[0,t_f]\to\R^3 \textrm{ measurable }\mid u(t)\in[0,1]^3 \textrm{ for almost every } t\}.
\end{equation*}
Values $0$ and $1$ for the controls correspond respectively to no feed or maximal feed rate.
We fix a positive threshold $\bar u\in(0,3t_f]$ and restrict the maximization problem \r{produtt} to controls $(u_1,u_2, u_3)\in{\cal U}$ which satisfy the further constraint 
 \bqnn
 \int_0^{t_f}u_1(s)+u_2(s)+u_3(s)\,ds=\bar u.
 \eqnn
The threshold $\bar u$ can be interpreted as the total amount of available substrates. Adding the constraint above means that we   exploit the totality of available substrates and oxygen. 

The problem shares a full coupling between a classical optimal control problem \r{eq:cs}, \r{produtt}, which describes the dynamics outside the cell, and a linear optimization problem \r{lp}, which models the metabolic activity inside the cell. In other words, at each time instant $t$, on the one hand, one needs to solve the linear program \r{lp} to obtain $\mu, v_e$ at time $t$ to plug into \r{eq:cs}; on the other hand, dynamics in   \r{eq:cs} must be integrated to in order to get the tuple $(G(t),Z(t),O(t),E(t))$ that allows to determine constraints in \r{lp}.

Since in   \r{lp} both the objective function and the constraints are linear with respect to $\bar v$, the outputs $\mu$ and $v_e$ are piecewise linear with respect to $\bar v$. Therefore  $\mu, v_e$ are  piecewise smooth  as functions of $(Z,G,O,E)$.
In the sequel we assume that  the outputs of \r{lp} are smooth. This amounts to  say  that the system evolves in a domain where $\mu, v_e$ are smooth as functions of $(Z,G,O,E)$. Under this assumption, we study optimal trajectories for  problem \r{eq:cs}, \r{produtt}. 

Let us give a more compact formulation of the optimization problem above. 
Let us rename 
\begin{equation}\label{sub}
 x_4=VG,~ x_5=VZ, ~x_6=VO, ~x_7=VE, ~x_8=VX,~ x_9=V.
\end{equation}
 In order   to get a cost of Mayer form in \r{produtt}, i.e., a cost depending only on the final point of the trajectory, it is natural to add three new variables to \r{eq:cs} by setting
\bqnn
x_i(t)=\int_0^tu_i(s)ds,~ i=1,2,3,
\eqnn 
so that the control system \r{eq:cs} reads
\begin{equation}
\begin{cases}
\dot x_1=u_1,\\
\dot x_2=u_2,\\
\dot x_3=u_3,\\
\dot x_4= F u_1-v_g(x_4,x_7,x_9)x_8,\\
\dot x_5= F u_2-v_z(x_4,x_5,x_7,x_9)x_8,\\
\dot x_6= F u_3-v_o(x_6,x_9)x_8,\\
\dot x_7=v_e(x_4,x_5,x_6,x_7,x_9)x_8,\\
\dot x_8=\mu(x_4,x_5,x_6,x_7,x_9)x_8,\\
\dot x_9=F,
\end{cases}\label{csnew}
\end{equation}
where
\begin{eqnarray}
v_g(x_4,x_7, x_9)&=&v_{g\mathrm{max}}\frac{x_4x_9}{(k_gx_9+x_4)(x_9+x_7/k^g_{i e})},\nonumber\\
v_z(x_4,x_5,x_7, x_9)&=&v_{z\mathrm{max}}\frac{x_5x_9^2}{(k_zx_9+x_5)(x_9+x_4/k_{i g})(x_9+x_7/k^z_{i e})},\label{eqve}\\
v_o(x_6,x_9)&=&v_{o\max}\frac{x_6}{k_ox_9+x_6},\nonumber
\end{eqnarray}
and $v_e(x_4,x_5,x_6,x_7,x_9), \mu(x_4,x_5,x_6,x_7,x_9)$ are defined in the obvious way by \r{lp}, \r{ve} using \r{sub}.

  The equation for $x_9$ can be trivially integrated. Indeed, the total volume of the culture plays the same role as time along the experiment. Here we prefer to keep the state variable $x_9$   so that the dynamics is  autonomous.

Set $x=(x_1,x_2,x_3,x_4,x_5,x_6, x_7,x_8,x_9)$, $u=(u_1,u_2,u_3)$, define the vector fields
\bqnn
F_0(x)=
\bpm 0\\0\\0\\-v_g(x_4,x_7,x_9)x_8\\-v_z(x_4,x_5,x_7,x_9)x_8\\-v_o(x_6,x_9)x_8\\v_e(x_4,x_5,x_6,x_7,x_9)x_8\\ \mu(x_4,x_5,x_6,x_7,x_9)x_8\\F
\epm, ~~
F_1(x)=
\bpm
1\\0\\0\\F\\0\\0\\0\\0\\0
\epm, ~~
F_2(x)=
\bpm
0\\1\\0\\0\\F\\0\\0\\0\\0
\epm,~~
F_3(x)=
\bpm
0\\0\\1\\0\\0\\F\\0\\0\\0
\epm,
\eqnn
the target
\bqnn
{\cal T}=\{x\in\R^9\mid x_1+x_2+x_3=\bar u\},
\eqnn 
and take an initial point 
\bqnn
x^0\in{\cal D}=\{x\in\R^9\mid x_1=x_2=x_3=0, x_9>0, x_i\geq 0, i=4 \dots 8\}.
\eqnn
 The 
   optimal control problem \r{produtt} (with the additional constraint on the controls) is equivalent to 
   \begin{equation}\label{ocp}
\begin{cases}
&\displaystyle{\max_{u\in{\cal U}}   \psi(x(t_f))} \\
&\dot x=f(x,u),\quad u\in{\cal U},\\
&x(0)=x^0,\quad x(t_f)\in {\cal T},
\end{cases}
\end{equation}
   where
$\psi(x) =x_7$ and
$
f(x,u)=F_0(x)+u_1F_1(x)+u_2F_2(x)
$.

We end this section by stating a result which ensures the existence of optimal solutions to \r{ocp} with suitable initial conditions.
\bp\label{p:ex}
Assume $v_e,\mu:\R^5\to\R$ are smooth, non negative and bounded from above.
  Then for every $x^0\in{\cal D}$, the optimal control problem \r{ocp} admits a solution. 
\ep
Note that, whenever $\mu, v_e$ are outputs of the optimization problem \r{lp}, they are non negative and bounded from above by construction, thanks to the constraints $0\leq \bar v_j\leq \tilde v_j$ in \r{lp}, so the assumptions of Proposition~\ref{p:ex} are satisfied in practical situations.

\medskip

\noindent{\it Proof.} Let $x^0\in{\cal D}$. We are going to prove that 
  every trajectory of the control system starting at $x^0$    is well-defined for all $t>0$ and satisfies 
  $$
  x(t)\in\tilde{\cal D}=\{x\in\R^9 \mid  x_i\geq 0, x_9>0\}, \quad\, t\geq 0.
  $$
Indeed we have $x_9(t)=x_9^0+F t\geq x_9^0>0$ for every $t$. 
Using the assumptions on $\mu, v_e$,  the dynamics is smooth on $\tilde {\cal D}$ and satisfies the sublinear growth condition 
\begin{equation}\label{sublinear}
|f(x,u)|\leq C(1+|x|),~~\forall\, (x,u)\in \tilde{\cal D}\times[0,1]^3,
\end{equation}
for a certain constant $C>0$. 
Let $T$ be such that a trajectory starting at $x^0$ is well-defined on $[0,T]$. Then $x(t)\in\tilde{\cal D}$ for every $t\in[0,T]$. Indeed, 
 inequalities $x_i(t)\geq0$ for $i=1,2,3$,  and $x_9(t)\geq x_9^0>0$ are trivially satisfied. Similarly,  since the equation for $x_8$ is linear with respect to $x_8$ and $
 \mu(x_4,x_5,x_6,x_7,x_9)\geq 0$,   we have  $x_8(t)\geq 0$. By assumption, the ethanol flux satisfies $v_e(x_4,x_5,x_6,x_7,x_9)\geq 0$,  whence $x_7(t)\geq 0$. We are left to prove that $x_i(t)\geq 0$ for $i=4,5,6$. By assumption $x_i(0)\geq 0$. Let $\bar t\geq 0$ be such that $x_i(\bar t)=0$. Then $\dot x_i(\bar t)= F u_i(\bar t)\geq 0$, hence $x_i(t)\geq 0$  for $t>\bar t$.  Eventually, trajectories starting at $x^0\in{\cal D}\subset\tilde {\cal D}$  satisfy $x(t)\in\tilde {\cal D}$ and, because of the growth condition \r{sublinear}, they  are well-defined for every $t>0$.

Choosing $u_1(s)\equiv u_2(s)\equiv u_3(s)\equiv \bar u /(3 t_f)$,    the trajectory associated to $(u_1,u_2, u_3)$ reaches the target ${\cal T}$ at time $t_f$. Hence for every $x^0\in {\cal D}$, a trajectory starting at $x^0$ and reaching the target always exists.
 The dynamics \r{csnew} is control-affine, whence the set of velocities $\{f(x,u)\mid u\in[0,1]^3\}$ is convex for every $x$. Moreover, the target ${\cal T}$ is closed, and the cost function $\psi$ is smooth. Therefore, the existence of an optimal solution is a consequence of classical results, see for instance \cite[Theorem~5.1.1]{benlibro}.
  \hfill$\blacksquare$

\section{Analysis of  extremal trajectories}\label{secmainth}

In this section we consider first order necessary conditions for optimality to identify properties of solutions to \r{ocp}. We first analyze a simplified case, where only one feeding substrate is consider and then turn to the case of three controls. 

Consider a general optimal control problem in Mayer form
\begin{equation}\label{ocpg}
\begin{cases}
&\displaystyle{\max_{u\in{\cal U}}   \psi(x(t_f))} \\
&\dot x=f(x,u),\quad u\in{\cal U},\\
&x(0)=x^0,\quad x(t_f)\in {\cal T},
\end{cases}
\end{equation}
   where  
$\psi\in\con^\infty(\R^N)$,  ${\cal U}=\{u:[0,t_f]\to[0,1]^m\mid u \textrm{ measurable}\}$, with control-affine dynamics
\begin{equation*}
f(x,u)=F_0(x)+\sum_{i=1}^mu_iF_i(x),
\end{equation*}
$F_i\in\con^\infty(\R^N\times\R^m,\R^N)$,
and ${\cal T}=\{x\mid \phi(x)=0\}$, with $\phi\in\con^\infty(\R^N)$.
The main   tool to look for solutions to \r{ocpg} is the Pontryagin Maximum Principle \cite{pontryagin-book} which is a first order necessary condition for optimality. 

The   Hamiltonian associated to \r{ocpg} is $H:\R^N\times\R^N\times\R^m\to\R$  
\bqnn
H(x, \lambda, u)=\la \lam, F_0(x)\ra+\sum_{i=1}^mu_i\la \lam, F_i(x)\ra,
\eqnn
where $\lam$ denotes the adjoint vector.
Let $u\in{\cal U}$ be an admissible control whose corresponding trajectory $x$ satisfies the terminal constraint $x(t_f)\in {\cal T}$.
The trajectory $x(\cdot)$ is called     {\it extremal} if there exist an absolutely continuous function  $\lam:[0,t_f]\to\R^N\setminus\{0\}$ and real numbers $(\nu_0,\nu_1)\neq (0,0)$, $\nu_0\geq 0$ such that $\lam(\cdot)$ satisfies the adjoint system
\begin{equation}\label{pmp1}
\dot\lam=-\frac{\partial H}{\partial x}(x,\lam,u), \quad \textrm{ for almost every }t\in[0,t_f],
\end{equation}
the maximisation condition
\begin{equation}\label{pmp2}
H{(x (t),\lam (t), u(t))}=\displaystyle{\max_{\omega\in[0,1]^m}\la \lam (t), F_0(x (t))\ra+\sum_{i=1}^m\omega_i\la \lam (t), F_i(x (t))\ra},~\textrm{ for almost every }t\in[0,t_f],
\end{equation}
and the transversality condition
\begin{equation}\label{pmp3}
\lam (t_f)=\nu_1\nabla\phi(x(t_f))+\nu_0\nabla\psi(x(t_f)).
\end{equation}
The  Pontryagin Maximum Principle, see for instance a version adapted to our problem in \cite[Theorem~6.3.1]{benlibro}, states that if $(u,x)$ is a solution to \r{ocpg}, then $x(\cdot)$ is an extremal trajectory.
Given an extremal trajectory,  define
 the switching functions
 \bqnn
 \ph_i(t)=\la \lam (t),F_i(x (t))\ra.
 \eqnn
Since the Hamiltonian is control-affine,  whenever $t$ is such that $\ph_i(t)\neq 0$, the maximization condition \r{pmp2} allows to identify $u_i(t)$ by $u_i(t)=\frac{1+\sign\ph_i(t)}{2}\in\{0,1\}$.  
We say that an extremal trajectory $x(\cdot)$ has a {\it singular arc} on $[T_0,T_1]\subset [0,t_f]$, with $T_0<T_1$, if the product $\Pi_{i=1}^m\ph_i(t)$ vanishes identically on  on $[T_0,T_1]$. Whenever $\ph_i(t)\equiv 0$ on $[T_0,T_1]$, we say that the control $u_i$ is {\it singular}.

We say that a measurable control $u_i:[0,t_f]\to[0,1]$ is {\it bang-bang} on $[T_0,T_1]$ if $u_i(t)\in\{0,1\}$ for almost every $t\in[T_0,T_1]$.  

The family of controls (and trajectories) among which one optimizes in  \r{ocpg} is infinite dimensional. Thus, in order to solve \r{ocpg} a fundamental step is to reduce the problem to a finite dimensional one. For instance, this can be done if, by means of necessary conditions, one deduces that every extremal trajectory is a finite concatenation of arcs where every control is smooth (e.g. bang-bang) and thus has a simple form. To do so, in the control-affine case, it suffices to prove that every extremal trajectory corresponds to a finite concatenation of bang-bang controls.
 In general this is hard to prove, mainly because either the concatenation of simple arcs may be infinite or because there may be singular arcs, where the control is not uniquely determined. 
 
 As concerns the first difficulty, it may indeed happen that an optimal trajectory features an accumulation of arcs where controls are constant. This is the well-known  Fuller  phenomenon, see \cite{fuller2}. In this work we do not deal with such pathological phenomena.  We refer the reader to \cite{fuller} for new ideas  overcoming numerical difficulties caused by Fuller phenomenon.

 As concerns the second issue, in the worst case, one can still expect that only a ``few'' extremal trajectories are ``bad'', i.e., they   feature a singular arc.   Roughly speaking, this amounts to prove that extremal trajectories corresponding to a singular control are  confined in some subsets of positive codimension. If so, one infers that optimal trajectories starting at most initial conditions do not have singular arcs and correspond to bang-bang controls.
 
 Withe regard to the presence of optimal singular trajectories for control-affine systems,   a genericity  result has been shown in \cite{cjt}. In that paper, the authors prove that, for a quadratic cost,  generically with respect to the control system there do not exist nontrivial optimal singular trajectories \cite[Corollary 2.9]{cjt}. Results of this type do  not apply to our problem because  system \r{csnew} is   not  generic,  $F_i$ being constant vector fields.  
  However, modeling $\mu, v_e$ as affine functions of $v_g, v_z, v_o$, in the three-input case, we are able to single out only two (out of seven) possible singular arcs and generically exclude all the other ones.

\subsection{Single-input case}\label{sinput}
In this section we study a simpler version of the optimal control problem  \r{ocp} by considering anaerobic growth and removing one feeding substrate, e.g., xylose. With this simplification, setting $y_2=VG,y_3=VE,y_4=VX,y_5=V$, the control system \r{csnew} becomes
\begin{equation}
\begin{cases}
\dot y_1=u\\
\dot y_2= F u-v_g(y_2,y_3,y_5)y_4\\
\dot y_3=v_e(y_2,y_3,y_5)y_4\\
\dot y_4=\mu(y_2,y_3,y_5)y_4\\
\dot y_5=F
\end{cases}\label{csuno}
\end{equation}
where $u\in {\cal U}=\{u:[0,t_f]\to\R \textrm{  measurable} \mid u(t)\in [0,1]\}$, 
\begin{equation}
v_g(y_2,y_3, y_5)=v_{g\mathrm{max}}\frac{y_2y_5}{(k_gy_5+y_2)(y_5+y_3/k^g_{i e})}\label{eqvguno}
\end{equation}
and $v_e(y_2,y_3,y_5), \mu(y_2,y_3, y_5)$ are given by
  \begin{eqnarray*}
&&\mu(y_2,y_3,y_5)=\max_{\bar v\in\R^{n}} \,\sum_{j=1}^{n} w_j \bar v_j\label{lp1}\\
&&\quad\qquad\qquad\qquad \mbox{ s.t. }  \, S\bar v = 0\nonumber\\
 &&\quad\qquad\qquad\qquad   \, \bar v_g=v_g(y_2,y_3,y_5) \nonumber\\
 &&\quad\qquad\qquad\qquad \, 0\leq  \bar v_j\leq \tilde v_j,\, j=1,\dots, n\nonumber\\
&& \nonumber\\
 && v_e(y_2,y_3,y_5) =\max \{\bar v_e\mid \bar v\in\mathrm{argmax }\, \mu(y_2,y_3,y_5)\}.\label{ve1}
 \end{eqnarray*}
System \r{csuno} is single-input control-affine, i.e., setting $y=(y_1,y_2,y_3,y_4,y_5)\in\R^5$,  the control system reads
\begin{equation}\label{csunfeed}
\dot y=g(y,u)=G_0(y)+u G_1(y), ~~y\in\R^5,
\end{equation}
where  
\begin{equation*}
G_0(y)=
\bpm 0\\-v_g(y_2,y_3,y_5)y_4\\v_e(y_2,y_3,y_5)y_4\\ \mu(y_2,y_3,y_5)y_4\\F
\epm, ~~
G_1(y)=
\bpm
1\\F\\0\\0\\0
\epm.
\end{equation*}
The target ${\cal T}$ becomes $\{y\in\R^5\mid \phi(y)=\bar u\}$, where $\bar u\in(0,t_f]$, $\phi(y)=y_1$, and the  function to be maximised is $\psi(y(t_f))$, 
with
$\psi(y) =y_3$. 
Therefore, given an initial condition $y^0\in{\cal D}=\{y\in\R^5\mid y_1^0=0, y_5^0>0, y_i\geq 0, i=2,3,4\}$, we consider the optimal control problem
\begin{equation}\label{ocp1}
\begin{cases}
&\displaystyle{\max_{u\in{\cal U}}   y_3(t_f)} \\
&\dot y=g(y,u),\quad u\in{\cal U},\\
&y(0)=y^0,\quad y(t_f)\in {\cal T}.
\end{cases}
\end{equation}
We model  the behavior of the parameters coming from the optimization problem according to the following assumption.
\begin{itemize}
\item[(H1)] There exist constants $a_1>0, b_1>0, \bar \mu\geq 0, \bar v\geq 0$ such that 
\begin{equation}\label{a1}
\begin{gathered}
\mu=a_1 v_g(y_2,y_3,y_5)+\bar \mu,\\
v_e=b_1 v_g(y_2,y_3,y_5)+\bar v.
\end{gathered}
\end{equation}
\end{itemize}
Under this hypothesis, $\mu$ and $v_e$ are smooth, non negative and bounded from above on the set $\{y\in\R^5\mid y_5>0, y_i\geq 0\}$, so that  
 existence of solutions to \r{ocp1}    can be proved as in  Proposition~\ref{p:ex}.  

The   Hamiltonian associated to \r{ocp1} is $H:\R^5\times\R^5\times\R\to\R$  
\bqnn
H(y, \eta, u)=\la \eta, G_0(y)\ra+u\la \eta, F(y)\ra=u(\eta_1+\eta_2 F)-\eta_2 v_g  y_4+\eta_3 v_e y_4+\eta_4 \mu  y_4+\eta F,
\eqnn
where $\eta$ denotes the adjoint vector.
Let $u\in{\cal U}$ be an admissible control such that $y$ is an extremal trajectory. Then  there exists an  absolutely continuous function  $\eta:[0,t_f]\to\R^5\setminus\{0\}$ which satisfies, for almost every $t\in[0,t_f]$, the adjoint system \r{pmp1}, i.e.,
\begin{equation}\label{syscovuno}
\begin{cases}
\dot\eta_1=0,\\
\dot \eta_2=y_4\left(\frac{\partial v_g}{\partial y_2}\eta_2-\frac{\partial v_e}{\partial y_2}\eta_3-\frac{\partial \mu}{\partial y_2}\eta_4\right),\\
\dot \eta_3=y_4\left(\frac{\partial v_g}{\partial y_3}\eta_2-\frac{\partial v_e}{\partial y_3}\eta_3-\frac{\partial \mu}{\partial y_3}\eta_4\right),\\
\dot \eta_4=v_g\eta_2-v_e\eta_3-\mu\eta_4,\\
\dot \eta_5=y_4\left(\frac{\partial v_g}{\partial y_5}\eta_2-\frac{\partial v_e}{\partial y_5}\eta_3-\frac{\partial \mu}{\partial y_5}\eta_4\right),
\end{cases}
\end{equation}
the maximisation condition \r{pmp2}
\begin{equation}\label{hh1}
H{(x (t),\eta (t), u(t))}=\displaystyle{\max_{0\leq \omega\leq 1}\la \eta (t), G_0(y (t))\ra+\omega\la \eta (t), G_1(y (t))\ra},
\end{equation}
and the transversality condition \r{pmp3}
\begin{equation}\label{finale}
\eta (t_f)=\nu_1\nabla\phi(y(t_f))+\nu_0\nabla\psi(y(t_f))=(\nu_1,0,\nu_0,0,0),
\end{equation}
for some constants $(\nu_0,\nu_1)\neq(0,0)$ with $\nu_0\geq 0$. 

 \begin{proposition}\label{unso}
Let $\mu, v_e$ satisfy assumption $(H1)$.  
Then any  optimal control for   \r{ocp1} is   bang-bang on $[0,t_f]$ and has no singular arcs.
 \end{proposition}
 
As a consequence, any solution $u(\cdot)$ to \r{ocp1} is given by a (possibly infinite) concatenation of arcs where $u(\cdot)$ is constantly equal to  $0$ or $1$.

\medskip

\noindent{\it Proof.} By the Pontryagin Maximum Principle, any optimal control for \r{ocp1} gives rise to an extremal trajectory.  Therefore, it suffices to show that given an extremal trajectory $y(\cdot)$ with adjoint vector  $\eta(\cdot)$ satisfying \r{syscovuno}, \r{hh1}, \r{finale}, then the switching function does not vanish identically on an interval of positive length.

The proof is split into two steps. First, we show that given an extremal trajectory $y(\cdot)$ then any solution $\eta(\cdot)$ of \r{syscovuno} which is constant on an interval $[T_0,T_1]\subset [0,t_f]$ with $T_1>T_0$ is constant on the whole interval $[0,t_f]$. Second, we prove that given an extremal trajectory $y(\cdot)$ and a corresponding covector satisfying \r{syscovuno}, \r{hh1}, \r{finale}, then the presence of a singular arc implies that the covector is constant on $[0,t_f]$ and we get a contradiction imposing the transversality condition.

\underline{Equilibria of \r{syscovuno}.} Let $u\in{\cal U}$ be an admissible control whose corresponding trajectory $y$ satisfies the terminal constraint $y(t_f)\in {\cal T}$ and is extremal. 
Let $\eta(\cdot)$ be a solution of \r{syscovuno}.
Once $y(\cdot)$ is given, 
system \r{syscovuno} can be written as $\dot\eta=M(y)\eta$ where $M(y)$ is
\begin{equation*}
M(y) = 
\begin{pmatrix}
    0 & 0 & 0 & 0 & 0 \\
    0 & y_4\frac{\partial v_g}{\partial y_2} & -b_1 y_4\frac{\partial v_g}{\partial y_2} & -a_1y_4\frac{\partial v_g}{\partial y_2} & 0 \\
0 & y_4\frac{\partial v_g}{\partial y_3} & -b_1 y_4\frac{\partial v_g}{\partial y_3} & -a_1y_4\frac{\partial v_g}{\partial y_3}& 0 \\
0 & v_g & -b_1 v_g-\bar v & -a_1v_g-\bar \mu & 0 \\
0 & y_4\frac{\partial v_g}{\partial y_5} & -b_1 y_4\frac{\partial v_g}{\partial y_5} & -a_1y_4\frac{\partial v_g}{\partial y_5}& 0 \\
\end{pmatrix}.
\end{equation*}
Reasoning as in the proof of Proposition~\ref{p:ex}, we infer that  $y(\cdot)$ is uniformly bounded on $[0,t_f]$ (growth in the control system \r{csuno} is  sublinear). As a consequence, the matrix function $t\mapsto M(y(t))$ is uniformly bounded on $[0,t_f]$. Hence we have existence and global uniqueness for solutions to \r{syscovuno}. 

Let us now analyze equilibria for \r{syscovuno}. To this aim, we need to compute the time dependent subspace $\ker M(y(t))$.
Let $t\in[0,t_f]$ and  $\bar \eta=(\bar \eta_1,\bar \eta_2,\bar \eta_3,\bar \eta_4,\bar \eta_5)\in \ker M(y(t))$. Then
$$
\begin{cases}
y_4\frac{\partial v_g}{\partial y_2}\left(\bar \eta_2-b_1\bar \eta_3-a_1\bar \eta_4\right)=0,\\
 y_4 \frac{\partial v_g}{\partial y_3}\left(\bar \eta_2-b_1\bar \eta_3-a_1\bar \eta_4\right)=0,\\
 v_g\left(\eta_2-b_1\bar \eta_3-a_1\bar \eta_4\right)-\bar v\bar \eta_3-\bar \mu\bar \eta_4=0,\\
y_4\frac{\partial v_g}{\partial y_5}\left(\bar \eta_2-b_1\bar \eta_3-a_1\bar \eta_4\right)=0.
\end{cases}
$$
Since $y_4(t)\frac{\partial v_g}{\partial y_2}(y_2(t),y_3(t),y_5(t))>0$ for every $t$, the system above is equivalent to
$$
\begin{cases}
\bar \eta_2-b_1\bar \eta_3-a_1\bar \eta_4=0,\\
\bar v\bar \eta_3+\bar \mu\bar \eta_4=0.
\end{cases}
$$
Therefore, $\ker M(y(t))$ does not depend on $t$ and coincides with
\begin{eqnarray}
\mathrm{Vect}\{(1,0,0,0,0),(0,b_1,1,0,0),(0,a_1,0,1,0),(0,0,0,0,1)\}, \textrm{ if } \bar\mu=\bar v=0 \label{kernel}\\ 
\mathrm{Vect}\{(1,0,0,0,0),(0,a_1\bar v-b_1\bar\mu,-\bar\mu,\bar v,0),(0,0,0,0,1)\}, \textrm{ otherwise.}\label{kernel2}
\end{eqnarray}
 Hence,  by global uniqueness, if a solution $\eta(\cdot)$ to \r{syscovuno} is constant on some interval $[T_0,T_1]$ with $T_1>T_0$ then $\eta(\cdot)$ is constant on the whole interval $[0, t_f]$ and belongs to the subspace defined in \r{kernel} or \r{kernel2}.

 \underline{Absence of singular arcs.}
 Let $u\in{\cal U}$ be an admissible control whose corresponding trajectory $y$ satisfies the terminal constraint $y(t_f)\in {\cal T}$ and is extremal.
Let $\eta(\cdot)$ be a solution of \r{syscovuno},\r{hh1},\r{finale}. Then $$\eta_1(t)\equiv \nu_1.$$
 By contradiction, assume that the switching function 
 $$
 \varphi(t)=\langle \eta(t), G_1(y(t))\rangle=\eta_1(t)+F\eta_2(t)
 $$ vanishes identically on an interval $[T_0, T_1]$. 
 Then $\eta_2$ is constant on $[T_0,T_1]$ and given by
\begin{equation}\label{qwer}
\eta_2(t)\equiv-\frac{\nu_1}{F}.
\end{equation}
Using the second equation in \r{syscovuno}, we obtain
\bqnn
F y_4\frac{\partial v_g}{\partial y_2}(\eta_2-b_1\eta_3-a_1\eta_4)\equiv 0
\eqnn
whence
\begin{equation}\label{daje1}
\eta_2(t)- b_1\eta_3(t) - a_1 \eta_4(t)\equiv 0.
\end{equation}
Differentiating \r{daje1} with respect to time, and taking account of \r{daje1}, we deduce 
\begin{equation}\label{rrrr}
a_1 (\eta_3(t) \bar v+\eta_4(t)\bar \mu )\equiv 0.
\end{equation}
Since $a_1>0$, conditions \r{daje1}, \r{rrrr} imply that $\eta(t)\in\ker M(y(t))$. Therefore,  $\eta(t)$ is constant on $[T_0,T_1]$   and henceforth on the whole interval $[0,t_f]$.
Moreover, by the transversality condition \r{finale}, the constant value of $\eta(t)$ is $\eta(t_f)=(\nu_1,0,\nu_0,0,0)$.  On one hand this implies that $\nu_1=0$, because of \r{qwer}. On the other hand, imposing  $(0,0,\nu_0,0,0)\in\ker M(y(t))$ implies that, either  $\nu_0 b_1 =0$ (when \r{kernel} applies) or 
 there exists $\alpha\in\R$ such that
$$
\begin{cases}
\alpha(a_1\bar v-b_1\bar \mu)=0,\\
-\alpha\bar\mu=\nu_0,\\
\alpha\bar v=0,
\end{cases}
$$
(when \r{kernel2} applies). In both cases,   we get $\nu_0=0$ which    contradicts $(\nu_1,\nu_0)\neq (0,0)$.
Hence, $\varphi$ does not vanish identically on a interval of positive length.
\hfill$\blacksquare$

\brem\label{remox} When we consider the single-input system obtained by removing oxygen and glucose, the conclusion of Proposition~\ref{unso} holds. This follows directly by the fact that, without oxygen and glucose, the equation for $v_z$ given in \r{eqve} becomes precisely \r{ve1} and the analysis is the same as the one above. Furthermore, when we remove glucose and xylose and keep oxygen as the only control, the analysis is even simpler, since $v_o$ in \r{eqve} does not depend on ethanol concentration. As a consequence, with the same technique, the conclusion of Proposition~\ref{unso} holds true for the single-input case where the only control is  oxygen feeding rate.
\erem

\subsection{Multi-input case}\label{minput}

 Let us go back to  the multi-input case where  different feeding substrates  are considered.
 The   Hamiltonian  associated to the optimal control problem \r{ocp} is $H:\R^9\times\R^9\times\R^3\to\R$ 
\bqnn
H(x, \lambda, u)=\la \lam, F_0(x)\ra+\sum_{i=1}^3u_i\la \lam, F_i(x)\ra,
\eqnn
where $\lambda\in\R^9$ denotes the adjoint vector.
We   model the behavior of  $\mu, v_e$ assuming the  following property.
\begin{itemize}
\item[(H2)] There exist constants $a_i, b_i>0$, $i=1,2,3$ and $\bar \mu, \bar v\geq 0$ such that 
\begin{equation}\label{a12}
\begin{gathered}
\mu=a_1 v_g+a_2 v_z+a_3 v_o+\bar \mu,\\
v_e=b_1v_g+b_2 v_z-b_3 v_o+\bar v.
\end{gathered}
\end{equation}
\end{itemize}
 Assuming $b_3>0$ is motivated by the fact that anaerobic conditions promote ethanol production, see \cite{biofuels}. Under this hypothesis, thanks to Proposition~\ref{p:ex} problem \r{ocp} admits a solution for every $x^0\in{\cal D}$.

 Let  $u\in{\cal U}$ be an admissible control such that the corresponding trajectory $x$ is extremal. Then
  there exists an  absolutely continuous function    $\lam:[0,t_f]\to\R^9\setminus\{0\}$ such that  
 \begin{eqnarray}
\dot\lam(t)&=&-\frac{\partial H}{\partial x}(x (t),\lam(t), u(t)), \textrm{ for almost every } t\in[0,t_f], \label{pmp}\\
\lam(t_f)&=&\nu_0\nabla\psi(x(t_f))+\nu_1\nabla\phi(x(t_f)),\label{trasv}
\end{eqnarray}
for some $(\nu_0,\nu_1)\neq(0,0)$ and $\nu_0\geq 0$, where $\psi(x)=x_7$ and $\phi(x)=x_1+x_2+x_3-\bar u$. Moreover,  for almost every $t\in[0,t_f]$, the triple $(x(t),\lam(t),u(t))$ satisfies  
\begin{equation}\label{hh}
H{(x (t),\lam (t), u(t))}=\max_{0\leq \omega_i\leq 1}\la \lam (t), F_0(x (t))\ra+\sum_{i=1}^3\omega_i\la \lam (t), F_i(x (t))\ra.
\end{equation}

\bt\label{mainthm}
Assume $\mu, v_e$ satisfy condition $(H2)$ and $a_1b_3-a_3b_1\neq 0$. 
Let  $(u(\cdot),x(\cdot))$ be a solution to \r{ocp}. Then, for every interval $[T_0,T_1]$,  at least one among  the controls  $u_1,u_2, u_3$ is bang-bang on $[T_0, T_1]$.
Moreover, if the trajectory has a singular arc on $[T_0,T_1]$, then one of the following conditions hold.
\bi
\iii[(i)] There exists a polyonomial $P:\R^9\to\R$ such that, for every $t\in [T_0,T_1]$, $x(t)\in\{x\in\R^9\mid P(x)=0\}$.
 \iii[(ii)] On $[T_0,T_1]$ the control $u_1$ is singular and the control $u_2$ is bang-bang.
\ei
 Finally,  under the additional assumption    
\begin{equation}\label{vzsemp}
v_z(x_5,x_7, x_9)=v_{z\mathrm{max}}\frac{x_5x_9}{(k_zx_9+x_5)(x_9+x_7/k^z_{i e})},
\end{equation}
only possibility (i) occurs.
\et
One possibility for a singular arc is that $u_1$ (which correspond to glucose feeding rate) is singular, $u_2$ is bang-bang ($u_3$ may or may not be singular). Otherwise, all other singular arcs are confined in some hypersurface of the type $\{x\in\R^9\mid P(x)=0\}$. Therefore, for most initial conditions, solutions to \r{ocp} correspond to controls $u$ such that $u_2$ is a concatenation of bang-bang arcs and $u_1, u_3$ are  concatenations of bang-bang or singular arcs.

 When we do not account for a preferred substrate and assume $v_z$ does not depend on glucose concentration $x_4$, then the result states that for most initial conditions optimal trajectories  for \r{ocp} correspond to concatenations of bang-bang arcs. In other words, generically with respect to initial conditions, singular arcs do not occur.

\medskip

\noindent{\it Proof of Theorem~\ref{mainthm}.}  
 By PMP  $x$ is an extremal trajectory, i.e., there exists an absolutely continuous $\lam:[0,t_f]\to\R^9\setminus\{0\}$  satisfying \r{pmp}, \r{trasv} for some $(\nu_0,\nu_1)\neq(0,0)$, $\nu_0\geq0$ and \r{hh}. 
The adjoint system in \r{pmp} reads
\begin{equation}\label{syscov}
\begin{cases}
\dot\lam_i=0,~~i=1,2,3\\
\dot \lam_4=x_8\left(\frac{\partial v_g}{\partial x_4}\lam_4+\frac{\partial v_z}{\partial x_4}\lam_5-\frac{\partial v_e}{\partial x_4}\lam_7-\frac{\partial \mu}{\partial x_4}\lam_8\right),\\
\dot \lam_5=x_8\left(\frac{\partial v_z}{\partial x_5}\lam_5-\frac{\partial v_e}{\partial x_5}\lam_7-\frac{\partial \mu}{\partial x_5}\lam_8\right),\\
\dot \lam_6=x_8\left(\frac{\partial v_o}{\partial x_6}\lam_6-\frac{\partial v_e}{\partial x_6}\lam_7-\frac{\partial \mu}{\partial x_6}\lam_8\right),\\
\dot \lam_7=x_8\left(\frac{\partial v_g}{\partial x_7}\lam_4+\frac{\partial v_z}{\partial x_7}\lam_5-\frac{\partial v_e}{\partial x_7}\lam_7-\frac{\partial \mu}{\partial x_7}\lam_8\right),\\
\dot \lam_8=v_g\lam_4+v_z\lam_5+v_o\lam_6-v_e\lam_7-\mu\lam_8,\\
\dot \lam_9=x_8\left(\frac{\partial v_g}{\partial x_9}\lam_4+\frac{\partial v_z}{\partial x_9}\lam_5+\frac{\partial v_o}{\partial x_9}\lam_6-\frac{\partial v_e}{\partial x_9}\lam_7-\frac{\partial \mu}{\partial x_9}\lam_8\right),
\end{cases}
\end{equation}
whereas the transversality condition \r{trasv} for the covector at final time takes the form  
\bqnn
\lam (t_f)=(\nu_1,\nu_1,\nu_1, 0,0,0, \nu_0,0,0).
\eqnn
The first three  equations in \r{syscov} are trivially integrated and give  
\begin{equation}
\lam_1(t)\equiv\lam_2(t)\equiv \lam_3(t)\equiv \nu_1.\label{eq:l1}
\end{equation}

The proof follows the idea of  Proposition~\ref{unso} and consists of 8 steps. Step 1 analyzes equilibria of \r{syscov} and ensures that
  any solution of \r{syscov} that is constant on some interval $[T_0,T_1]\subset [0,t_f]$ is actually constant on the whole interval $[0,t_f]$. Step 2 is concerned with singular arcs where all the switching functions vanish identically and implies that at least one among the controls is bang-bang. Steps 3,4,5 deal with singular arcs where two switching functions vanish identically and one does not. Steps 6,7,8 deal with singular arcs where one switching function vanishes identically and the other two  do not.

    \smallskip

 \underline{Step 1.} Equilibria of \r{syscov} satisfy
\begin{equation}
\begin{cases}
x_8\left(\frac{\partial v_g}{\partial x_4}(\lam_4-b_1\lam_7-a_1\lam_8)+\frac{\partial v_z}{\partial x_4}(\lam_5-b_2\lam_7-a_2\lam_8)\right)=0,\\
x_8\frac{\partial v_z}{\partial x_5}(\lam_5-b_2\lam_7-a_2\lam_8)=0,\\
x_8\frac{\partial v_o}{\partial x_6}(\lam_6+b_3\lam_7-a_3\lam_8)=0,\\
x_8\left(\frac{\partial v_g}{\partial x_7}(\lam_4-b_1\lam_7-a_1\lam_8)+\frac{\partial v_z}{\partial x_7}(\lam_5-b_2\lam_7-a_2\lam_8)\right)=0,\\
v_g(\lam_4-b_1\lam_7-a_1\lam_8)+v_z(\lam_5-b_2\lam_7-a_2\lam_8)+v_o(\lam_6+b_3\lam_7-a_3\lam_8)-(\bar v\lam_7+\bar\mu\lam_8)=0,\\
x_8\left(\frac{\partial v_g}{\partial x_9}(\lam_4-b_1\lam_7-a_1\lam_8)+\frac{\partial v_z}{\partial x_9}(\lam_5-b_2\lam_7-a_2\lam_8)+\frac{\partial v_o}{\partial x_9}(\lam_6+b_3\lam_7-a_3\lam_8)\right)=0.
\end{cases}
\end{equation}
By the second equation, since $x_8\frac{\partial v_z}{\partial x_5}>0$, we deduce $\lam_5-b_2\lam_7-a_2\lam_8=0$. Plugging this information into the first equation, since $x_8\frac{\partial v_g}{\partial x_4}>0$, we infer that $\lam_4-b_1\lam_7-a_1\lam_8=0$. By the third equation, since $x_8\frac{\partial v_o}{\partial x_6}>0$, we obtain $\lam_6+b_3\lam_7-a_3\lam_8=0$. Then system above is equivalent to
\begin{equation}
\begin{cases}
\lam_4-b_1\lam_7-a_1\lam_8=0,\\
\lam_5-b_2\lam_7-a_2\lam_8=0,\\
\lam_6+b_3\lam_7-a_3\lam_8=0,\\
\bar v\lam_7+\bar\mu\lam_8=0.
\end{cases}
\end{equation}
Therefore, a vector $\lambda\in\R^9$ is an equilibrium if and only if 
\begin{eqnarray} 
\lambda\in\mathrm{Vect}\{e_1,e_2,e_3,e_9,(0,0,0,b_1,b_2,-b_3,1,0,0),(0,0,0,a_1,a_2,a_3,0,1,0)\}, \textrm{ if } \bar \mu=\bar v= 0,\label{cher}\\
\lambda\in\mathrm{Vect}\{e_1,e_2,e_3,e_9,(0,0,0,b_1\bar\mu-a_1\bar v,b_2\bar\mu-a_2\bar v, -b_3\bar\mu-a_3\bar v, \bar\mu, -\bar v, 0)\}, \textrm{ otherwise,}\label{cher2}
\end{eqnarray}
where $e_1=(1,0,0,0,0,0,0,0,0)$ and $e_2,e_3,e_9$ are defined accordingly. 
In the sequel we set
$
A=\lam_4-b_1\lam_7-a_1\lam_8,
B=\lam_5-b_2\lam_7-a_2\lam_8,
C=\lam_6+b_3\lam_7-a_3\lam_8,
D=\bar v\lam_7+\bar\mu\lam_8$.
As it happens, $\lam$ is en equilibrium if and only if $A=B=C=D=0$, and \r{syscov} writes
\begin{equation}\label{syscovabc}
\begin{cases}
\dot\lam_i=0, ~~i=1,2,3\\
\dot \lam_4=x_8\left(\frac{\partial v_g}{\partial x_4} A+\frac{\partial v_z}{\partial x_4}B\right),\\
\dot \lam_5=x_8\frac{\partial v_z}{\partial x_5}B,\\
\dot \lam_6=x_8 \frac{\partial v_o}{\partial x_6}C,\\
\dot \lam_7=x_8\left(\frac{\partial v_g}{\partial x_7}A+\frac{\partial v_z}{\partial x_7}B\right),\\
\dot \lam_8=v_gA+v_zB+v_oC-D,\\
\dot \lam_9=x_8\left(\frac{\partial v_g}{\partial x_9}A+\frac{\partial v_z}{\partial x_9}B+\frac{\partial v_o}{\partial x_9}C\right).
\end{cases}
\end{equation}

\smallskip

\underline{Step 2.}
Let $[T_0,T_1]\subset [0,t_f]$ be a singular arc such that $\varphi_i|_{[T_0,T_1]}\equiv 0$, $i=1,2,3$. Then, $\lam_i$ is constant for $i=4,5,6$ and given by $\lam_i(t)=-\frac{\nu_1}{F}$.
 Imposing that  the right-hand sides of the  fourth, fifth and sixth equations in \r{syscovabc} vanish identically we deduce that $A=B=C=0$.
Plugging this information into the  the seventh and ninth equations in \r{syscovabc}  gives that $\lam_7,\lam_9$ are also constant on $[T_0,T_1]$. As a consequence, by condition $\lam_4-b_1\lam_7-a_1\lam_8=0$ we obtain that $\lam_8$ is constant on $[T_0,T_1]$. Imposing  that the right-hand side of the equation for $\lam_8$ vanishes we conclude that $D=0$.
Therefore, for every $t\in[T_0,T_1]$, $\lam(t)$ belongs to the subspace defined in \r{cher} or \r{cher2}. By step 1, $\lam(t)$ is constant on $[T_0,T_1]$ and  henceforth on the whole interval $[0,t_f]$.
Imposing the transversality condition, one gets $\lam(t)\equiv \lam(t_f)$, which immediately gives $\nu_1=0$. Finally,  $\lam(t)=(0,0,0,0,0,0,\nu_0,0,0)$ is an equilibrium of \r{syscov} if and only if
 either $\nu_0 b_i=0$ (when \r{cher} applies) or   there exists $\alpha\in\R$ such that 
$$
\begin{cases}
\alpha(b_1\bar \mu-a_1\bar v)=0,\\
\alpha(b_2\bar \mu-a_2\bar v)=0,\\
\alpha(b_3\bar \mu+a_3\bar v)=0,\\
\alpha\bar\mu=\nu_0,\\
\alpha\bar v=0,
\end{cases}
$$
(when \r{cher2} applies). In both cases  we immediately get $\nu_0=0$, which gives a contradiction with $(\nu_1,\nu_0)\neq (0,0)$. As a consequence, we conclude that on any subinterval $[T_0,T_1]$, at least one control amongst $u_1, u_2,u_3$ is bang-bang.

\smallskip

\underline{Step 3.}
Let now $[T_0,T_1]\subset[0,t_f]$ be a singular arc where $\varphi_1$ is nonzero and $\varphi_2,\varphi_3$ vanish identically. Then $u_1(t)=\frac{1+\sign\varphi_1(t)}{2}\in\{0,1\}$ for almost every $t$.
Moreover, for every $t\in[T_0,T_1]$
\bqnn
\lam_5(t)\equiv\lam_6(t)\equiv -\frac{\nu_1}{F}. 
\eqnn
Conditions $\dot\varphi_2\equiv 0$, $\dot \varphi_3\equiv 0$, which are equivalent to $\dot\lam_5\equiv 0$,  $\dot\lam_6\equiv 0$, respectively, imply 
$
B=C=0$.  Imposing that $\frac{dB}{dt}\equiv \frac{dC}{dt}\equiv 0$ and taking account of  \r{syscovabc} and of $B=C=0$, gives rise to the following linear system in $A, D$  
\begin{equation}\label{asdf}
\begin{cases}
\left(b_2 x_8\frac{\partial v_g}{\partial x_7}+a_2 v_g\right)A-a_2 D=0,\\
\left(b_3 x_8\frac{\partial v_g}{\partial x_7}-a_3 v_g\right) A+a_3 D=0.
\end{cases}
\end{equation}
Up to a positive constant, the determinant of system above is $\frac{\partial v_g}{\partial x_7}$ which vanishes only when $x_4(t)=0$. Define $P(x)=x_4$.
 Assume there exists a $t\in[T_0,T_1]$ such that $P(x(t))\neq 0$. Then the determinant of system above is nonzero, whence $A=D=0$. By step 1, $\lam$ is an equilibrium of \r{syscov} and therefore is constant on $[0,t_f]$. Using the transversality condition,  $\lam(t)\equiv \lam(t_f)=(\nu_1,\nu_1,\nu_1,0,0,0,\nu_0,0,0)$. From the fifth component, we deduce  immediately that $\nu_1=0$. Imposing that $(0,0,0,0,0,0,\nu_0,0,0)$ is an equilibrium, either $\nu_0 b_i=0$ 
 (when \r{cher} applies), or  there exists $\g\in\R$ such that
\begin{equation*}
\begin{cases}
\g(a_2\bar v-b_2\bar\mu)=0,\\
\g(a_1\bar v-b_1\bar\mu)=0,\\
-\g\bar\mu=\nu_0,\\
\g\bar v=0,
\end{cases}
\end{equation*}
(when \r{cher2} applies).
As in step 2, in both cases we obtain $\nu_0$ and thus a contradiction with $(\nu_0,\nu_1)\neq (0,0)$. 

Assume now that $P(x(t))=x_4(t)\equiv 0$ on $[T_0, T_1]$. Since $x_4(\cdot)$ satisfies $\dot x_4 = F u_1-v_g x_8$, due to the form of $v_g$ as a function of $x_4$ (see \r{eqve}), the only possibility for this to happen is that  $x_4^0=0$ and $u_1(t)=0$ on $[0,T_1]$ and in this case   the trajectory is contained in $\{x\in\R^9\mid P(x)=0\}$.
Moreover,   $\lam_i|_{[T_0,T_1]}$ is constant for every $i\neq 4$. Indeed, for $x_4\equiv 0$,
both equations in \r{asdf} are equivalent to $D=0$. Using \r{syscovabc},  we obtain that $\lam_7,\lam_8,\lam_9$ are constant on $[T_0,T_1]$ and the equation for $\lam_4(t)$ depends only on $x(t)$ and on the constant values of $\lam_8,\lam_9$.

\smallskip

\underline{Step 4.} 
Let   $[T_0,T_1]\subset[0,t_f]$ be a singular arc where $\varphi_3$ is nonzero and $\varphi_1,\varphi_2$ vanish identically. Then $u_3(t)=\frac{1+\sign\varphi_3(t)}{2}\in\{0,1\}$ for almost every $t$.
Moreover, 
\bqnn
\lam_4(t)\equiv\lam_5(t)\equiv -\frac{\nu_1}{F}, t\in[T_0,T_1]. 
\eqnn
Imposing that the right-hand sides of equations for $\lam_4,\lam_5$ in \r{syscovabc} vanish identically we get $A=B=0$. Then, the right-hand side of the equation for $\lam_7$ also vanishes identically, that is, $\lam_7$ is constant on $[T_0,T_1]$. By condition $A=0$, this implies that $\lam_8$ is constant on $[T_0,T_1]$, whence 
\begin{equation}\label{cidi}
v_oC-D=0.
\end{equation}
 Differentiating this condition with respect to time we deduce 
$$
\left(\frac{d}{dt}v_o + x_8 v_o \frac{\partial v_o}{\partial x_6}\right) C=0.
$$
By the form of $v_o$ (see \r{eqve}), $\frac{d}{dt}v_o$ depends only on $x_6, x_9$ and $\dot x_6,\dot x_9$. Now, $\dot x_9=F$ and $\dot x_6= Fu_3-v_o x_8$. Since $u_3$ is constant on $[T_0,T_1]$, the coefficient 
$$
\frac{d}{dt}v_o + x_8 \frac{\partial v_o}{\partial x_6}=\frac{Fv_{o\max}k_o}{(k_ox_9+x_6)^2}(u_3 x_9-x_6),
$$ is a rational function of $(x_6,x_9)$ and does not depend explicitly on $t$.  
Define $P(x)=u_3x_9-x_6$. Assume there exists $t\in[T_0,T_1]$ such that $P(x(t))\neq 0$. Then $C=0$ and, by \r{cidi}, $D=0$. Therefore $\lam$ is an equilibrium of \r{syscov}. By step 1, $\lam$ is constant on $[0,t_f]$  and belongs to the subspace defined either in \r{cher} or \r{cher2}. Reasoning as in step 3, we get a contradiction imposing the transversality condition. Otherwise,  for every $t\in[T_0,T_1]$, $x(t)\in\{x\in\R^9\mid P(x)=0\}$. In this case, there are two possibilities : $u_3\equiv 1$ or $u_3\equiv 0$. If $u_3\equiv 1$ then $P(x(t))\equiv 0$ implies $x_6(t)\equiv x_9(t)$. Recalling that $x_6,x_9$ satisfy
$$
\begin{cases}
\dot x_6=F u_3-v_{o\max}\frac{x_6x_8}{k_ox_9+x_6},\\
\dot x_9=F,
\end{cases}
$$
we have a contradiction ($x_8(t)>0$ for every $t\geq T_0$). Hence the only possibility is that $u_3\equiv 0$. Then $P(x(t))\equiv 0$ implies $x_6(t)\equiv 0$ on $[T_0,T_1]$. Because of the evolution equation for $x_6$ this means that $x_6^0=0$ and $u_3(t)\equiv 0$ for every $t\in[0,T_1]$. 

\smallskip

\underline{Step 5.} Let   $[T_0,T_1]\subset[0,t_f]$ be a singular arc where $\varphi_2$ is nonzero and $\varphi_1,\varphi_3$ vanish identically. Then $u_2(t)=\frac{1+\sign\varphi_2(t)}{2}\in\{0,1\}$ for almost every $t$.
Moreover, 
\bqnn
\lam_4(t)\equiv\lam_6(t)\equiv -\frac{\nu_1}{F}, t\in[T_0,T_1]. 
\eqnn
Imposing that $\lam_4,\lam_6$ are constant, by \r{syscovabc} we deduce $C=0$ and
\begin{equation}\label{pabd}
\frac{\partial v_g}{\partial x_4}A+\frac{\partial v_z}{\partial x_4}B=0.
\end{equation}
Differentiating condition $C=0$ with respect to time and taking \r{pabd}, \r{syscovabc} into account we obtain another equation in $(A,B,D)$, namely
\begin{equation}\label{dueabd}
\left(b_3x_8\frac{\partial v_g}{\partial x_7}-a_3 v_g\right)A+\left(b_3x_8\frac{\partial v_z}{\partial x_7}-a_3 v_z\right)B+a_3 D=0
\end{equation}
In order to carry out the same argument as in steps 3 and 4, one should find a third linear and homogeneous condition in $A,B,D$. The problem is that, if one differentiates \r{pabd} or \r{dueabd} with respect to time, the obtained conditions involve $\dot x_4$, which depends on the control $u_1$. Since we are on an arc where $\varphi_1\equiv 0$, this control is singular and thus unknown. Therefore, singular arcs of this type cannot be excluded in general.  

On the contrary, we can do better under the additional assumption that 
 $v_z$ is given as in \r{vzsemp}. In this case, since $v_z$ does not depend in $x_4$, condition \r{pabd} is equivalent to $A=0$. Differentiating conditions $A=0, C=0$ with respect to time we obtain a linear homogeneous system in $B,D$, namely
\begin{equation}\label{sssss}
\begin{cases}
\left(b_1 x_8 \frac{\partial v_z}{\partial x_7}+ a_1 v_z\right) B-a_1 D=0,\\
\left(b_3 x_8 \frac{\partial v_z}{\partial x_7}+ a_3 v_z\right) B-a_3 D=0.
\end{cases}
\end{equation}
The determinant of the  system above is
$(a_1b_3-a_3b_1)x_8\frac{\partial v_z}{\partial x_7},
$
only depends on $(x_4,x_7,x_8,x_9)$ and does not depend explicitly on $t$. Since $a_1b_3-a_3b_1\neq 0$, the determinant is nonzero if and only if $x_5(t)=0$.   If there exists $t\in[T_0,T_1]$ such that $x_5(t)\neq 0$ then   $\lam$ is constant on $[0,t_f]$ and this leads to a contradiction with the transversality condition. Assume now that $x_5\equiv 0$ on $[T_0,T_1]$. Because of the equation for $x_5$, this can only happen if   $u_2\equiv 0$ on $[0,T_1]$.
 In this case the trajectory is contained in $\{x\in\R^9\mid P(x)=0\}$ with $P(x)=x_5$. Moreover, reasoning as at the end of  step 3, $\lam_i$ is constant on $[T_0,T_1]$ for every $i\neq 5$ and the equation for $\lam_5$ only depends on $x(t)$ and the constant values of $\lam_7,\lam_8$.


\smallskip

\underline{Step 6.} Let   $[T_0,T_1]\subset[0,t_f]$ be a singular arc where $\varphi_1,\varphi_3$ are nonzero and $\varphi_2$ vanishes identically. Then $u_i(t)=\frac{1+\sign\varphi_i(t)}{2}\in\{0,1\}$, $i=1,3$ for almost every $t$.
Moreover, 
\bqnn
\lam_5(t)\equiv -\frac{\nu_1}{F}, t\in[T_0,T_1]. 
\eqnn
From \r{syscovabc} we infer that $B=0$. Differentiating this condition with respect to time we obtain
\begin{equation}\label{uuna}
\left(b_2 x_8\frac{\partial v_g}{\partial x_7}+a_2 v_g\right) A+ a_2 v_o C -a_2 D =0.
\end{equation}
Differentiating twice \r{uuna} with respect to time gives rise to other two linear homogeneous equations in $(A,C,D)$ whose coefficients depend only on $(x_4,x_6,x_7,x_8,x_9)$ and $(\dot x_4,\dot x_6,\dot x_7,\dot x_8,\dot x_9)$. Since the controls $u_1, u_3$ are constant on $[T_0,T_1]$, these coefficients do not depend explicitly on $t$. Considering the obtained linear homogeneous system in $(A,C,D)$ there are two possibilities. Either there exists $t\in [T_0, T_1]$ such that the determinant of the system does not vanish or the determinant vanishes identically. In the first case we deduce that $\lam$ is an equilibrium of \r{syscov} and we get a contradiction with the transversality condition. In the second case, there exists a polynomial function $P:\R^9\to \R$ depending only on $(x_4,x_6,x_7,x_8,x_9)$ such that the trajectory is contained in $\{x\in\R^9\mid P(x)=0\}$.

\smallskip

\underline{Step 7.} Let   $[T_0,T_1]\subset[0,t_f]$ be a singular arc where $\varphi_1,\varphi_2$ are nonzero and $\varphi_3$ vanishes identically. Then $u_i(t)=\frac{1+\sign\varphi_i(t)}{2}\in\{0,1\}$, $i=1,2$ for almost every $t$.
Moreover, 
\bqnn
\lam_6(t)\equiv -\frac{\nu_1}{F}, t\in[T_0,T_1]. 
\eqnn
From \r{syscovabc} we obtain that $C=0$. Differentiating this condition with respect to time we obtain
\begin{equation}\label{uuuna}
\left(b_3 x_8\frac{\partial v_g}{\partial x_7}-a_3 v_g\right) A+ \left(b_3 x_8\frac{\partial v_z}{\partial x_7}-a_3 v_z\right) B +a_3 D =0.
\end{equation}
Differentiating twice \r{uuuna} with respect to time gives rise to other two linear homogeneous equations in $(A,B,D)$ whose coefficients depend only on $(x_4,x_5,x_7,x_8,x_9)$ and $(\dot x_4,\dot x_5,\dot x_7,\dot x_8,\dot x_9)$. Since the controls $u_1, u_2$ are constant on $[T_0,T_1]$, these coefficients do not depend explicitly on $t$. Thus, as in step 6,  we conclude that either such a singular arc does not occur or   there exists a polynomial function $P:\R^9\to \R$ depending only on $(x_4,x_5,x_7,x_8,x_9)$ such that the trajectory is contained in $\{x\in\R^9\mid P(x)=0\}$.

\smallskip

\underline{Step 8.} Let   $[T_0,T_1]\subset[0,t_f]$ be a singular arc where $\varphi_2,\varphi_3$ are nonzero and $\varphi_1$ vanishes identically. Then $u_i(t)=\frac{1+\sign\varphi_i(t)}{2}\in\{0,1\}$, $i=2,3$ for almost every $t$.
Moreover, 
\bqnn
\lam_4(t)\equiv -\frac{\nu_1}{F}, t\in[T_0,T_1]. 
\eqnn
Using \r{syscovabc} this implies
\begin{equation}\label{trio}
\frac{\partial v_g}{\partial x_4}A+\frac{\partial v_z}{\partial x_4}B=0.
\end{equation}
As is the case in step 5, differentiating \r{trio} with respect to time provides  linear homogeneous conditions in $(A, B,C,D)$ whose coefficients involve $\dot x_4= F u_1 -v_g x_8$. Since the control $u_1$ is singular and thus unknown, we are not able to exclude singular arcs of this type in general.

Under the additional assumption that $v_z$ is as in \r{vzsemp}, the situation simplifies. More precisely, $\lam_4$ constant implies $A=0$. The system 
$$
\begin{cases}
\frac{dA}{dt}=0,\\
\frac{d^2A}{dt^2}=0,\\
\frac{d^3A}{dt^3}=0,
\end{cases}
$$
is linear and homogeneous in $(B,C,D)$ and its coefficients depend only on $(x_5,x_6,x_7,x_8,x_9)$ and 
$(\dot x_5,\dot x_6,\dot x_7,\dot x_8,\dot x_9)$. Recalling that the controls $u_2, u_3$ are constant on $[T_0,T_1]$ we can conclude as in step 7 and prove that either such a singular arc does not occur or the trajectory is contained in a hypersurface $\{x\in\R^9\mid P(x)=0\}$, where $P$ is a polynomial function depending only on $(x_5,x_6,x_7,x_8,x_9)$.
  \hfill$\blacksquare$

\brem
Concerning singular extremals having 2 singular controls and a bang-bang control there are three possibilities, that we consider   in steps 3, 4, 5 of the proof above. 
In the first and second cases (steps 3, 4 respectively), the only possibility for a singular extremal to exist is that glucose, respectively oxygen, concentration is identically zero since the beginning of the experiment. In view of applications, these situations are rather meaningless, being that glucose is the preferred substrate for biomass growth and aerobic conditions enhance biomass growth.   In the last case (step 5) under the simplified assumption \r{vzsemp}, the only possibility is that xylose concentration is identically zero since the beginning of the experiment. This represents the only reasonable possibility from the point of view of applications.   
 \erem

 \section{Conclusions}\label{conclu}
  
In \cite{biofuels}, the authors perform two series of {\it in silico} experiments for {\it Saccharomyces cerevisiae}. The first one considers glucose media, whereas the second one deals with glucose and xylose media. In both series, feeding rates (which in our framework are represented by $u_1,u_2$) are assumed constant and ethanol productivity is maximized among trajectories characterized by a dissolved oxygen concentration of $50\%$ on $[0,t_s]$ and of $0\%$ from $[t_s,t_f]$ and the switching time $t_s$ is the only control parameter. In Section~\ref{sinput} we consider the more general case where the feeding rate of a substrate is treated as a control. More precisely, the analysis carried out   for the single-input case
shows that when the oxygen feeding rate is treated as control then optimal trajectories are bang-bang (see Remark~\ref{remox}). This provides the theoretical background needed to  legitimate the choice in \cite{biofuels} of piecewise constant dissolved oxygen concentration: optimizing ethanol productivity among all possible input profiles for the oxygen feeding rate is equivalent to optimizing among feeding rates that are bang-bang. Furthermore, in Section~\ref{minput} we consider a unified model where all feeding rates (glucose, xylose, oxygen) are treated as controls and $\mu, v_e$ are affine functions of $v_g,v_z,v_o$. In this more general framework, we prove that, at least when there is no preferred substrate (among glucose and xylose), then  optimal trajectories exiting from most initial conditions are characterized bang-bang
feeding rates.  In particular, in this simplified context, the only case of interest (in view of applications) for a singular extremal having two singular controls is when glucose and oxygen are singular and xylose concentration is zero since the beginning of the experiment. It would be of interest to check optimality of these particular singular extremals via {\it in silico} experiments.

Finally, the analysis of singular extremals provided here constitutes the building block for the development of the hybrid model where parameters $\mu$ and $v_e$ are only piecewise smooth. To go further, the starting point is to study necessary conditions for optimality, which, for this generalized context, can be found in \cite{garavello-piccoli-hyb,sussmann-hybrid}.

\small{
\bibliographystyle{abbrv}
\bibliography{biblio_biofuels}

\begin{thebibliography}{10}

\bibitem{Alford}
J.~Alford.
\newblock Bioprocess control: Advances and challenges.
\newblock {\em Computers \& Chemical Engineering}, 30(10-12):1464--1475, 2006.

\bibitem{optcontbiofi}
P.~T. Benavides and U.~Diwekar.
\newblock Optimal control of biodiesel production in a batch reactor: Part i:
  Deterministic control.
\newblock {\em Fuel}, 94(0):211 -- 217, 2012.

\bibitem{hyb}
M.~S. Branicky.
\newblock Introduction to hybrid systems.
\newblock In {\em Handbook of networked and embedded control systems}, Control
  Eng., pages 91--116. Birkh\"auser Boston, Boston, MA, 2005.

\bibitem{benlibro}
A.~Bressan and B.~Piccoli.
\newblock {\em Introduction to the mathematical theory of control}, volume~2 of
  {\em AIMS Series on Applied Mathematics}.
\newblock American Institute of Mathematical Sciences (AIMS), Springfield, MO,
  2007.

\bibitem{gauthier-busvelle}
{\'E}.~Busvelle and J.-P. Gauthier.
\newblock On determining unknown functions in differential systems, with an
  application to biological reactors.
\newblock {\em ESAIM Control Optim. Calc. Var.}, 9:509--551, 2003.

\bibitem{fuller}
M.~Caponigro, R.~Ghezzi, B.~Piccoli, and E.~Tr\'elat.
\newblock {R}egularization of chattering phenomena via bounded variation
  control.
\newblock preprint 2013, arXiv:1303.5796.

\bibitem{cjt}
Y.~Chitour, F.~Jean, and E.~Tr{\'e}lat.
\newblock Singular trajectories of control-affine systems.
\newblock {\em SIAM J. Control Optim.}, 47(2):1078--1095, 2008.

\bibitem{rfba}
M.~W. Covert, C.~Schilling, and B.~Palsson.
\newblock Regulation of gene expression in flux balance models of metabolism.
\newblock {\em Biophys J.}, 83(3):1331--1340, 2002.

\bibitem{ifba}
M.~W. Covert, N.~Xiao, T.~J. Chen, and J.~R. Karr.
\newblock Integrating metabolic, transcriptional regulatory and signal
  transduction models in {\it {e}scherichia {c}oli}.
\newblock {\em Bioinformatics}, 24(18):2044--2050, 2008.

\bibitem{hyblibro}
M.~D. Di~Benedetto and A.~Sangiovanni-Vincentelli.
\newblock {\em Hybrid Systems: Computation and Control}.
\newblock Lecture Notes in Comput. Sci. 2034. Springer-Verlag, Berlin,
  Heidelberg, 2001.

\bibitem{optcontbiof}
T.~Eevera, K.~Rajendran, and S.~Saradha.
\newblock Biodiesel production process optimization and characterization to
  assess the suitability of the product for varied environmental conditions.
\newblock {\em Renewable Energy}, 34(3):762 -- 765, 2009.

\bibitem{fuller2}
A.~T. Fuller.
\newblock Study of an optimum non-linear control system.
\newblock {\em J. Electronics Control (1)}, 15:63--71, 1963.

\bibitem{garavello-piccoli-hyb}
M.~Garavello and B.~Piccoli.
\newblock Hybrid necessary principle.
\newblock {\em SIAM J. Control Optim.}, 43(5):1867--1887 (electronic), 2005.

\bibitem{gauthier-bior}
J.-P. Gauthier, H.~Hammouri, and S.~Othman.
\newblock A simple observer for nonlinear systems applications to bioreactors.
\newblock {\em IEEE Trans. Automat. Control}, 37(6):875--880, 1992.

\bibitem{hensoneth}
J.~L. Hjersted and M.~A. Henson.
\newblock Optimization of fed-batch saccharomyces cerevisiae fermentation using
  dynamic flux balance models.
\newblock {\em Biotechnol. Prog.}, 22:1239--1248, 2006.

\bibitem{hensonethanol}
J.~L. Hjersted and M.~A. Henson.
\newblock Steady-state and dynamic flux balance analysis of ethanol production
  by saccharomyces cerevisiae.
\newblock {\em IET Systems Biology}, 3:167--179, 2009.

\bibitem{biofuels}
J.~L. Hjersted, M.~A. Henson, and R.~Mahadevan.
\newblock Genome-{S}cale {A}nalysis of {\it {s}accharomyces cervisiae}
  {M}etabolism and {E}thanol {P}roduction in {F}ed-{B}atch {C}ulture.
\newblock {\em Biotechnology and Bioengineering}, 97(5):1190--1204, 2007.

\bibitem{tb}
E.~Jung, S.~Lenhart, and Z.~Feng.
\newblock Optimal control of treatments in a two-strain tubercolosis model.
\newblock {\em Discrete and Continuous Dynamical Systems--Series B},
  2(4):473--482, November 2002.

\bibitem{hiv}
D.~Kirschner, S.~Lenhart, and S.~Serbin.
\newblock Optimal control of the chemotherapy of {HIV}.
\newblock {\em J. Math. Biol.}, 35:775--792, 1997.

\bibitem{odemodel}
A.~Kremling, K.~Bettenbrock, and E.~Gilles.
\newblock Analysis of global control of {E}scherichia coli carbohydrate uptake.
\newblock {\em BMC Systems Biology}, 1(42), 2007.

\bibitem{dfba}
R.~Mahadevan, J.~Edwards, and F.~r. Doyle.
\newblock Dynamic flux balance analysis of diauxic growth in {E}scherichia
  coli.
\newblock {\em J Theor Biol.}, 213(1):73--88, 2001.

\bibitem{moreno}
J.~Moreno.
\newblock Optimal time control of bioreactors for the wastewater treatment.
\newblock {\em Optimal Control Applications Methods}, 20(3):145--164, 1999.

\bibitem{palsson-book}
B.~O. Palsson.
\newblock {\em Systems Biology - Property of Reconstructed Networks}.
\newblock Cambridge University Press, 2006.

\bibitem{pontryagin-book}
L.~S. Pontryagin, V.~G. Boltyanski{\u\i}, R.~V. Gamkrelidze, and E.~F.
  Mishchenko.
\newblock {\em The {M}athematical {T}heory of {O}ptimal {P}rocesses}.
\newblock ``Nauka'', Moscow, fourth edition, 1983.

\bibitem{Rapaport}
A.~Rapaport and D.~Dochain.
\newblock Minimal time control of fed-batch processes with growth functions
  having several maxima.
\newblock {\em IEEE Trans. Automat. Contr.}, 56(11):2671--2676, 2011.

\bibitem{sussmann-hybrid}
H.~J. Sussmann.
\newblock A nonsmooth hybrid maximum principle.
\newblock In {\em Stability and stabilization of nonlinear systems ({G}hent,
  1999)}, volume 246 of {\em Lecture Notes in Control and Inform. Sci.}, pages
  325--354. Springer, London, 1999.

\bibitem{vaccine}
S.~Tiwari, P.~Verma, P.~Singh, and R.~Tuli.
\newblock Plants as bioreactors for the production of vaccine antigens.
\newblock {\em Biotechnology Advances}, 27(4):449--467, 2009.

\bibitem{bioprocess-review}
K.~Yamuna~Rani and V.~S. Ramachandra~Rao.
\newblock Control of fermenters - a review.
\newblock {\em Bioprocess and Biosystems Engineering}, 21:77--88, 1999.
\newblock 10.1007/PL00009066.

\end{thebibliography}
}

\end{document}